\documentclass[12pt]{article}

\usepackage{amsmath}
\usepackage{amsfonts}
\usepackage{amsthm}
\usepackage{hyperref}
\usepackage{enumitem}
\usepackage{array}
\usepackage[nottoc,numbib]{tocbibind}
\usepackage{xr}
\newcommand{\dimh}{\mbox{$\dim_{\mathrm{H}}$}}
\newcommand{\dm}{\mbox{$\mu^{\mathrm{H}}$}}

\newcommand{\R}{\mbox{$\mathbb{R}$}}
\newcommand{\N}{\mbox{$\mathbb{N}$}}

\begin{document}


\author{Attila Losonczi}
\title{The generalized Hausdorff-integral and its applications}

\date{\today}

\newtheorem{thm}{\qquad Theorem}[section]
\newtheorem{prp}[thm]{\qquad Proposition}
\newtheorem{lem}[thm]{\qquad Lemma}
\newtheorem{cor}[thm]{\qquad Corollary}
\newtheorem{rem}[thm]{\qquad Remark}
\newtheorem{ex}[thm]{\qquad Example}
\newtheorem{df}[thm]{\qquad Definition}
\newtheorem{prb}{\qquad Problem}

\newtheorem*{thm2}{\qquad Theorem}
\newtheorem*{cor2}{\qquad Corollary}
\newtheorem*{prp2}{\qquad Proposition}

\maketitle

\begin{abstract}

\noindent

We go on to widen the scope of the previously defined Hausdorff-integral in the sense that we extend it to functions taking values in $[0,+\infty)\times[0,+\infty)$. In all our intentions, we follow the same attitude that we had in our previous investigations, i.e. we work in the realm of Hausdorff dimension and measure.

\noindent
\footnotetext{\noindent
AMS (2020) Subject Classifications:  28A25, 28A78\\

Key Words and Phrases: Hausdorff dimension and measure, generalized integral in measure spaces, generalized measures}

\end{abstract}


\section{Introduction}

This paper can be considered as a continuation of the papers \cite{lahi} and \cite{lahi2}. For original intentions, aims and motivation, we suggest that the reader should consult to \cite{lahi} and \cite{lahi2}.

In \cite{lahi} we defined the Hausdorff-integral on $\mathbb{R}$ for real functions, then in \cite{lahi2} we generalized it for real functions on so-called h-measure spaces and found connection between the two integral notions.
We can go on generalizing and define the integral of functions whose values are taken from $[0,+\infty)\times[0,+\infty)$ and are defined on h-measure spaces. Our applications will be based on the integral of such functions. 

The new integral has the usual properties, e.g. it is a linear functional (homogeneous and additive), monotone, it provides an h-measure, etc.
However this generalization of the integral has some flaws. It is not true anymore that if $f_n\to f,\ (0,0)\leq f_n\leq f_{n+1}\ (n\in\N)$, then $(\mathrm{H})\hspace{-0.4pc}\int_Kf_n\ d\mu\to(\mathrm{H})\hspace{-0.4pc}\int_Kf\ d\mu$. It is a consequence of the fact that the topological space $[0,+\infty)\times[0,+\infty]$ is not separable. It is not true either that every measurable function is a limit of simple functions, therefore in some of the proofs, different, more sophisticated methods have to be applied.

\smallskip

The structure of the paper is the following. First, we define a multiplication in $[0,+\infty)\times[-\infty,+\infty]$, as it is essential for the integral. Afterwards, we define the new integral for simple functions (some call them step functions) and analyze its properties. In the next step, we extend the integral for non-negative measurable functions and again we derive its properties. In many of the cases, our proofs follow the usual arguments however we always have to listen to the slight but in some cases important differences. Finally, we present some applications which are sophisticated deficiency measurements.

\subsection{Basic notions and notations}\label{sbnn}

Here we enumerate the basics that we will apply throughout the paper.

\smallskip

For $K\subset\mathbb{R}$, $\chi_K$ will denote the characteristic function of $K$ i.e. $\chi_K(x)=1$ if $x\in K$, otherwise $\chi_K(x)=0$.

\smallskip

We will use the usual notations $\mathbb{N}_0=\mathbb{N}\cup\{0\},\ \overline{\mathbb{N}}_0=\mathbb{N}\cup\{0,+\infty\},\newline \mathbb{R}^+=\{x\in\mathbb{R}:x>0\},\ \mathbb{R}^+_0=\mathbb{R}^+\cup\{0\},\ \overline{\mathbb{R}}^+_0=\mathbb{R}^+_0\cup\{+\infty\},\ \overline{\mathbb{R}}=\mathbb{R}\cup\{-\infty,+\infty\}.$

\smallskip

Some usual operations and relation with $\pm\infty$: $(+\infty)+(+\infty)=+\infty,\ (-\infty)+(-\infty)=-\infty$; if $r\in\mathbb{R}$, then $r+(+\infty)=+\infty,\ r+(-\infty)=-\infty,\ {-\infty<r<+\infty}$. $+\infty+(-\infty)$ is undefined. If $r>0$, then $r\cdot (+\infty)=+\infty,\ r\cdot(-\infty)=-\infty$. If $r<0$, then $r\cdot (+\infty)=-\infty,\ r\cdot(-\infty)=+\infty$. And $0\cdot\infty=0$.

\smallskip

Set $S(x,\delta)=\{y \in\mathbb{R}^n:d(x,y)<\delta\},\ \ \overset{.}{S}(x,\delta)=S(x,\delta)-\{x\}\ \ (x\in\mathbb{R}^n,\delta>0)$.

\smallskip

In $\mathbb{R}^n\ \pi_k$ denotes the \textbf{projection} to the $k^{th}$ coordinate, in other words $\pi_k(x_1,\dots,x_n)=x_k\ \ (1\leq k\leq n)$.

\smallskip

If $K\subset\mathbb{R}^n$, then $\dim_{\mathrm{H}}(K)$ will denote the \textbf{Hausdorff dimension} of $K$, and $\mu^{d}(K)$ will denote the $d$ dimensional \textbf{Hausdorff measure} of $K$.

\smallskip

$\lambda$ will denote the Lebesgue-measure and also the outer measure as well.

\smallskip
 
 We call $f:K\to\R$ a \textbf{simple function}, if $\mathrm{Ran }f$ is finite.

\subsection{Basic definitions from \cite{lahi} and \cite{lahi2}}

For easier readability, we copy some basic definitions from \cite{lahi} and \cite{lahi2} here. 

\smallskip

If $(d_1,m_1),(d_2,m_2),\in[0,+\infty)\times[-\infty,+\infty]$, then let
\[\boldsymbol{(d_1,m_1)\leq(d_2,m_2)}\iff d_1<d_2\text{ or }(d_1=d_2\text{ and }m_1\leq m_2).\]

\smallskip

In the sequel, we will apply the order topology on $[0,+\infty)\times[-\infty,+\infty]$.

\smallskip

If $(d_1,m_1),(d_2,m_2),\in[0,+\infty)\times[-\infty,+\infty]$, then let
\[\boldsymbol{(d_1,m_1)+(d_2,m_2)}=\begin{cases}
(d_2,m_2)&\text{if }d_1<d_2\\
(d_1,m_1)&\text{if }d_2<d_1\\
(d_1,m_1+m_2)&\text{if }d_1=d_2\text{ and }m_1+m_2\text{ exists}\\
\text{undefined}&\text{if }d_1=d_2\text{ and }m_1+m_2\text{ is undefined}.
\end{cases}\]

\smallskip

If $(d,m)\in[0,+\infty)\times[-\infty,+\infty]$ and $c\in\R$, then let 
\[\boldsymbol{c(d,m)}=\begin{cases}
(d,cm)&\text{if }c\ne 0\\
(0,0)&\text{if }c=0.
\end{cases}\]

\smallskip

We define $\dm:\mathrm{P}(\mathbb{R})\to\{0\}\times\overline{\mathbb{N}}_0\cup(0,1]\times[0,+\infty]$. If $K\subset\mathbb{R}$, then $\boldsymbol{\dm(K)}=(d,m)\text{ where }d=\dimh K,m=\mu^d(K).$

\smallskip

\label{dHint}Let $K\subset\mathbb{R}$ be Borel-measurable and $f:K\to\mathbb{R}$ be a Borel-measurable function. Set $D_f=\{x\in K:f(x)\ne 0\}$. Set $\boldsymbol{(\mathrm{H})\hspace{-0.4pc}\int_Kf}=(d,m)\in[0,1]\times[-\infty,+\infty],\text{\ where\ \ }d=\dimh D_f,\ m=\int_{D_f}f\ d\mu^d,$
and we use the conventions that $\dimh\emptyset=0$ and $\int_{\emptyset}f\ d\mu^d=0$.

\smallskip

Let ${\cal S}\subset\mathrm{P}(K)$ be a $\sigma$-algebra. We call $\mu:{\cal S}\to[0,+\infty)\times[0,+\infty]$ an \textbf{h-measure} on ${\cal S}$, if $\mu(\emptyset)=(0,0)$, and $(A_n)$ being a sequence of sets from ${\cal S}$ such that $A_n\cap A_m=\emptyset\ (n\ne m)$ implies that $\mu\left(\bigcup_{n=1}^{\infty}A_n\right)=\sum_{n=1}^{\infty}\mu(A_n).$
In this case we call $(K,{\cal S},\mu)$ an \textbf{h-measure space}.

\smallskip

In the sequel, in this subsection, $(K,{\cal S},\mu)$ will denote an h-measure space. Let $f:K\to\R$ be a measurable non-negative simple function, i.e. $f=\sum_{i=1}^nc_i\chi_{A_i}$ where $c_i\in\R^+_0,\ A_i\in{\cal S}$  $(1\leq i\leq n)$ and $A_i\cap A_j=\emptyset\ (i\ne j)$. Then set $\boldsymbol{(\mathrm{H})\hspace{-0.4pc}\int_Kf\ d\mu}=\sum_{i=1}^nc_i\mu(A_i).$

\smallskip

Let $f:K\to\R$ be a measurable non-negative function. Then set $\boldsymbol{(\mathrm{H})\hspace{-0.4pc}\int_Kf\ d\mu}=\sup\left\{(\mathrm{H})\hspace{-0.4pc}\int_Kg\ d\mu:g\in\Phi_f\right\}$, where \newline $\Phi_f=\{g:0\leq g\leq f,\ g\text{ is a measurable simple function}\}$.

\smallskip

Let $f:K\to\R$ be a measurable non-positive function. Then set $\boldsymbol{(\mathrm{H})\hspace{-0.4pc}\int_Kf\ d\mu}=(-1)\cdot(\mathrm{H})\hspace{-0.4pc}\int_K-f\ d\mu.$

\smallskip

Let $f:K\to\R$ be a measurable function. Let $f^+=\max\{f,0\}, f^-=\min\{f,0\}$. Then set $\boldsymbol{(\mathrm{H})\hspace{-0.4pc}\int_Kf\ d\mu}=(\mathrm{H})\hspace{-0.4pc}\int_Kf^+\ d\mu  +  (\mathrm{H})\hspace{-0.4pc}\int_Kf^-\ d\mu$, if the sum exists.

\section{Generalized integral}

\subsection{Multiplication on $[0,+\infty)\times[-\infty,+\infty]$}

In order to define the integral, we will need multiplication on $[0,+\infty)\times[0,+\infty)$. We will use a natural way to do that, however one might argue against it, saying that that rule does not always hold for Hausdorff dimension and measure. We admit that it is true, but we could not find a better way until now, and more importantly, this method of multiplication and the resulting integral suit our later purposes completely.

\smallskip

We define the operation on a greater domain that we will actually need.

\begin{df}Let $(d_1,m_1),(d_2,m_2)\in[0,+\infty)\times[-\infty,+\infty]$. Let 
\[(d_1,m_1)(d_2,m_2)=\begin{cases}
(0,0)&\text{if }(d_1,m_1)=(0,0)\text{ or }(d_2,m_2)=(0,0)\\
(d_1+d_2,m_1m_2)&\text{otherwise.}\tag*{\qed}
\end{cases}\]
\end{df}

\begin{rem}\label{r2}The multiplication of a real number and a value from $[0,+\infty)\times[-\infty,+\infty]$ is a special case of the multiplication defined here, because for ${c\in\R}, (d,m)\in[0,+\infty)\times[-\infty,+\infty]$, we have $c(d,m)=(0,c)(d,m)=(d,cm)$ if $c\ne 0$, while if $c=0$, then both sides equal to $(0,0)$. Therefore if we identify $c$ with $(0,c)$, then we can embed the former algebraic structure into the current one.
\end{rem}

\begin{prp}\label{p14}Let $(d_1,m_1),(d_2,m_2)\in[0,+\infty)\times[-\infty,+\infty]$. Then $(d_1,m_1)(d_2,m_2)=(0,0)$ iff $(d_1,m_1)=(0,0)$ or $(d_2,m_2)=(0,0)$.
\end{prp}
\begin{proof}The sufficiency is obvious.
\par To see the necessity, suppose that none of terms is $(0,0)$. Then we get that $d_1+d_2=0$ which gives that $d_1=0$ and $d_2=0$. Furthermore $m_1m_2=0$ yields that $m_1=0$ or $m_2=0$. Hence one of the terms is $(0,0)$ which is a contradiction.
\end{proof}

\begin{prp}The multiplication is commutative and associative.
\end{prp}
\begin{proof}Being commutative is trivial.
\par For the associative property, let $(d_1,m_1),(d_2,m_2),(d_3,m_3)\in[0,+\infty)\times[-\infty,+\infty]$. If any of the three terms is $(0,0)$, then clearly the result is $(0,0)$ in both cases.
If $(d_1,m_1)\cdot(d_2,m_2)=(0,0)$, then $(d_1,m_1)=(0,0)$ or $(d_2,m_2)=(0,0)$; similarly for $(d_2,m_2)\cdot(d_3,m_3)$. If none of the three terms is $(0,0)$, then $(d_1,m_1)\big((d_2,m_2)(d_3,m_3)\big)=(d_1+d_2+d_3,m_1m_2m_3)=\big((d_1,m_1)(d_2,m_2)\big)(d_3,m_3)$.
\end{proof}

\begin{prp}\label{p9}The distributive law holds in $[0,+\infty)\times[0,+\infty]$.
\end{prp}
\begin{proof}Let $(d_1,m_1),(d_2,m_2),(d_3,m_3)\in[0,+\infty)\times[0,+\infty]$. We have to show that $(d_1,m_1)\big((d_2,m_2)+(d_3,m_3)\big)=(d_1,m_1)(d_2,m_2)+(d_1,m_1)(d_3,m_3)$. There are four cases to be checked.
\vspace{-0.6pc}\begin{enumerate}\setlength\itemsep{-0.3em}
\item If $(d_1,m_1)=(0,0)$, then we are done.
\item If either $(d_2,m_2)=(0,0)$ or $(d_3,m_3)=(0,0)$, then it holds again. 
\item If $(d_2,m_2)+(d_3,m_3)=(0,0)$, then $d_2=d_3=0$ and $m_2=m_3=0$. 
\item If none of the above hold, then let $d=\max\{d_2,d_3\}$. Then we get that
\[(d_1,m_1)\big((d_2,m_2)+(d_3,m_3)\big)=\left(d_1+d;m_1\sum_{\substack{i=2,3 \\ d_i=d}}m_i\right),\]
\[(d_1,m_1)(d_2,m_2)+(d_1,m_1)(d_3,m_3)=(d_1+d_2,m_1m_2)+(d_1+d_3,m_1m_3)=\]
\[\left(d';\sum_{\substack{i=2,3 \\ d_1+d_i=d'}}m_1m_i\right)=\left(d_1+d;m_1\sum_{\substack{i=2,3 \\ d_i=d}}m_i\right),\]
where $d'=\max\{d_1+d_2,d_1+d_3\}$, and using that $d'=d_1+d$.\qedhere
\end{enumerate}
\end{proof}

\begin{ex}The distributive law does not hold in general if we allow negative numbers in the second coordinate as well. 
\par See e.g. $(1,1)\big((0,5)+(0,-5)\big)=(0,0)\ne(1,5)+(1,-5)=(1,0)$.\qed
\end{ex}

\begin{prp}\label{p10}Let a sequence $\big((d_i,m_i)\big)$ be given on $[0,+\infty)\times[0,+\infty]$ and let $(d,m)\in[0,+\infty)\times[0,+\infty]$. Then
\[(d,m)\sum_{i=1}^{\infty}(d_i,m_i)=\sum_{i=1}^{\infty}(d,m)(d_i,m_i).\]
\end{prp}
\begin{proof}If $(d,m)=(0,0)$, then we are done. If $\sum_{i=1}^{\infty}(d_i,m_i)=(0,0)$, then all terms equal to $(0,0)$, and we are done again. Also note, that the terms equal to $(0,0)$ can be left out without changing the results of both sides. If none of the previous conditions hold, then let $D=\sup\{d_i:i\in\mathbb{N}\}$. By \cite{lahi} 2.12,
\[(d,m)\sum_{i=1}^{\infty}(d_i,m_i)=(d,m)\left(D;\sum_{d_j=D}m_j\right)=\left(d+D;m\sum_{d_j=D}m_j\right),\]
while
\[\sum_{i=1}^{\infty}(d,m)(d_i,m_i)=\sum_{i=1}^{\infty}(d+d_i,mm_i)=\left(D';\sum_{d+d_j=D'}mm_j\right)=\left(d+D;\sum_{d_j=D}mm_j\right),\]
where $D'=\sup\{d+d_i:i\in\mathbb{N}\}$, and using that $D'=d+D$.
\end{proof}

\begin{rem}We note that the distributive laws (\ref{p9} and \ref{p10}) hold in a more general context, namely the multiplier ($(d_1,m_1)$ and $(d,m)$) can be from $[0,+\infty)\times[-\infty,+\infty]$, but the terms inside the sum must be from the reduced set $[0,+\infty)\times[0,+\infty]$.
\end{rem}

\begin{prp}\label{p11}Let $(d_1,m_1),(d_2,m_2),(d_3,m_3)\in[0,+\infty)\times[0,+\infty]$ such that $(d_1,m_1)\leq(d_2,m_2)$. Then $(d_1,m_1)(d_3,m_3)\leq(d_2,m_2)(d_3,m_3)$.
\end{prp}
\begin{proof}From the condition we get that $d_1<d_2\text{ or }(d_1=d_2\text{ and }m_1\leq m_2)$. If either $(d_3,m_3)=(0,0)$ or $(d_1,m_1)=(0,0)$, then we are done. Similarly when $(d_2,m_2)=(0,0)$, because that implies that $(d_1,m_1)=(0,0)$. Hence we can suppose that none of the terms is $(0,0)$. Therefore one has to check that $(d_1+d_3,m_1m_3)\leq(d_2+d_3,m_2m_3)$ which holds in both cases mentioned in the first sentence of the proof.
\end{proof}

\subsection{Measurability}

First we need to define measurability of generalized functions.

\begin{df}Let $(K,{\cal S})$ be a measurable space (i.e. ${\cal S}$ is a $\sigma$-algebra on $\mathrm{P}(K)$). Let ${f:K\to[0,+\infty)\times[0,+\infty)}$ be a function. We say that $f$ is \textbf{measurable} if for each $(d,m)\in[0,+\infty)\times[0,+\infty)$ $\{x\in K:f(x)<(d,m)\}\in{\cal S}$ holds.\qed
\end{df}

\begin{prp}Let $(K,{\cal S})$ be a measurable space and ${f:K\to[0,+\infty)\times[0,+\infty)}$. Then the following statements are equivalent.
\vspace{-0.6pc}\begin{enumerate}\setlength\itemsep{-0.3em}
\item $f$ is measurable.
\item For each $(d,m)\in[0,+\infty)\times[0,+\infty]$, $\{x\in K:f(x)\leq(d,m)\}\in{\cal S}$.
\item For each $(d,m)\in[0,+\infty)\times[0,+\infty]$, $\{x\in K:(d,m)<f(x)\}\in{\cal S}$.
\item For each $(d,m)\in[0,+\infty)\times[0,+\infty]$, $\{x\in K:(d,m)\leq f(x)\}\in{\cal S}$.
\end{enumerate}
\end{prp}
\begin{proof}All equivalences simply follow from the fact that the space $[0,+\infty)\times[0,+\infty]$ is first countable and T2.
\end{proof}

\begin{prp}Let $(K,{\cal S})$ be a measurable space and ${f:K\to[0,+\infty)\times[0,+\infty)}$ be measurable. Then the following statements hold.
\vspace{-0.6pc}\begin{enumerate}\setlength\itemsep{-0.3em}
\item $\pi_1(f)$ is measurable.
\item For all $d\in[0,+\infty)$ $\pi_2(f|_{K_d})$ is measurable, where $K_d=\{x\in K:\pi_1(f(x))=d\}$.
\end{enumerate}
\end{prp}
\begin{proof}\begin{enumerate}\setlength\itemsep{-0.3em}
\item Clearly $\{x:\pi_1(f(x))<d\}=\{x:f(x)<(d,0)\}\ \ (d\in[0,+\infty))$.
\item Obviously $\{x:\pi_2(f|_{K_d})<m\}=\{x:\pi_1(f(x))=d,\pi_2(f(x))<m\}=\{x:(d,0)\leq f(x)<(d,m)\}\ \ (m\in[0,+\infty))$.\qedhere
\end{enumerate}
\end{proof}

\begin{ex}If $f$ is measurable, then $\pi_2(f)$ is not necessarily measurable.
\par Let $K=[0,1], {\cal S}$ be the Borel sets, $g:[0,1]\to[0,1]$ be a non Borel measurable function and $f(x)=(x,g(x))$ when $x\in[0,1]$. Clearly $\{x\in K:f(x)<(d,m)\}$ equals to either $\{x\in K:x<d\}$ or $\{x\in K:x\leq d\}$, and both sets are Borel, hence $f$ is measurable.
\end{ex}

\begin{prp}Let $(K,{\cal S})$ be a measurable space and ${f:K\to[0,+\infty)\times[0,+\infty)}$ be measurable and $(d',m')\in[0,+\infty)\times[0,+\infty)$. Then $(d',m')f$ is measurable as well.
\end{prp}
\begin{proof}Let $(d,m)\in[0,+\infty)\times[0,+\infty)$.
\par If $(d',m')=(0,0)$, then the statement is trivial.
\par If $m'>0$, then clearly $\{x:(d',m')f(x)\leq(d,m)\}={\{x:f(x)\leq(d-d',\frac{m}{m'})\}}$, because if $f(x)=(d'',m'')$, then $(d',m')f(x)=(d'+d'',m'm'')\leq(d,m)$ holds iff $d''<d-d'$ or $d''=d-d'$ and $m'm''\leq m$ iff $(d'',m'')\leq(d-d',\frac{m}{m'})$.
\par If $m'=0,d'>0$, then $\{x:(d',m')f(x)\leq(d,m)\}={\{x:f(x)\leq(d-d',+\infty)\}}$, where similar argument works.
\end{proof}

\begin{prp}Let $(K,{\cal S})$ be a measurable space and ${f,g:K\to[0,+\infty)\times[0,+\infty)}$ be measurable. Then $f+g$ is measurable as well.
\end{prp}
\begin{proof}Let $(d,m)\in[0,+\infty)\times[0,+\infty)$. Then we simply get that
\[\{x\in K:f(x)+g(x)<(d,m)\}=\]
\[\big(\{x:f(x)<(d,m)\}\cap\{x\in K:g(x)<(d,0)\}\big)\bigcup\big(\{x:f(x)<(d,0)\}\cap\{x:g(x)<(d,m)\}\big)\bigcup\]
\[\{x:(d,0)\leq f(x)<(d,m),\ (d,0)\leq g(x)<(d,m),\ f(x)+g(x)<(d,m)\}.\]
The last set can be written as
\[\{x:(d,0)\leq f(x)<(d,m),\ (d,0)\leq g(x)<(d,m),\ \pi_2(f(x))+\pi_2(g(x))<m\}=\]
\[\bigcup_{q_r}\{x:(d,0)\leq f(x)<(d,m),\ (d,0)\leq g(x)<(d,m),\ \pi_2(f(x))<q_r<m-\pi_2(g(x))\}=\]
\[\bigcup_{q_r}\big(\{x:(d,0)\leq f(x)<(d,q_r)\}\cap\{x:(d,0)\leq g(x)<(d,m-q_r)\}\big),\]
where $(q_r)$ is an enumeration of $\mathbb{Q}$.
\end{proof}

\begin{prp}Let $(K,{\cal S})$ be a measurable space and ${f_n:K\to[0,+\infty)\times[0,+\infty)}$ be measurable ($n\in\N$), and let $f_n\to f$. Then $f$ is measurable as well.
\end{prp}
\begin{proof}Let $(d_1,m_1);(d_2,m_2)\in[0,+\infty)\times[0,+\infty)$, $(d_1,m_1)<(d_2,m_2)$. For any function $g$, and for $n\in\N$ let 
\[N^n_g=\left\{x\in K:\left(d_1-\frac{1}{n},0\right)\leq g(x)\leq\left(d_2,m_2+\frac{1}{n}\right)\right\},\]
if $m_1=0$, and let
\[N^n_g=\left\{x\in K:\left(d_1,m_1-\frac{1}{n}\right)\leq g(x)\leq\left(d_2,m_2+\frac{1}{n}\right)\right\},\]
if $m_1>0$.
It can be readily checked that 
\[\{x\in K:(d_1,m_1)\leq f(x)\leq(d_2,m_2)\}=\bigcap_{n=1}^{\infty}\bigcup_{m=1}^{\infty}\bigcap_{k=m}^{\infty}N^n_{f_k}.\qedhere\]
\end{proof}

\subsection{The integral of functions taking values in ${[0,+\infty)\times[0,+\infty)}$}

We are now ready to define the integral in the promised more general context.

\begin{df}\label{d5}Let $(K,{\cal S},\mu)$ be an h-measure space. Let ${f:K\to[0,+\infty)\times[0,+\infty)}$ be a measurable simple function, i.e. $f=\sum_{i=1}^n(d_i,m_i)\chi_{A_i}$, where $(d_i,m_i)\in[0,+\infty)\times[0,+\infty),\ A_i\in{\cal S}$  $(1\leq i\leq n)$ and ${A_i\cap A_j=\emptyset}\ {(i\ne j)}$. Then set
\[(\mathrm{H})\hspace{-0.4pc}\int\limits_Kf\ d\mu=\sum_{i=1}^n(d_i,m_i)\mu(A_i).\tag*{\qed}\]
\end{df}

When building the theory of the integral, we then can follow the same way we had in subsection \cite{lahi2}\ref{hi2-ss1} and we can show that the integral of simple functions is well-defined, and having the required properties (see \cite{lahi2}\ref{hi2-p6}), and the integral provides an h-measure (see \cite{lahi2}\ref{hi2-p3}). For those, one has to observe only that in this more general context too, we only need that the multiplication is commutative, associative and distributive, and the properties in \ref{p14} and \ref{p11} hold. We just repeat those three statements without proofs; the reader can check them easily.

\begin{prp}The integral in definition \ref{d5} is well defined.\qed
\end{prp}

\begin{prp}\label{p6b}Let $(K,{\cal S},\mu)$ be an h-measure space. Let $f,g:K\to[0,+\infty)\times[0,+\infty)$ be measurable simple functions. Then the followings hold.
\vspace{-0.6pc}\begin{enumerate}\setlength\itemsep{-0.3em}
\item $(0,0)\leq(\mathrm{H})\hspace{-0.4pc}\int_Kf\ d\mu$.
\item If $f\leq g$, then $(\mathrm{H})\hspace{-0.4pc}\int_Kf\ d\mu\leq(\mathrm{H})\hspace{-0.4pc}\int_Kg\ d\mu$.
\item If $c\in[0,+\infty)\times[0,+\infty)$, then $(\mathrm{H})\hspace{-0.4pc}\int_Kcf\ d\mu=c\cdot(\mathrm{H})\hspace{-0.4pc}\int_Kf\ d\mu$.
\item $(\mathrm{H})\hspace{-0.4pc}\int_Kf+g\ d\mu=(\mathrm{H})\hspace{-0.4pc}\int_Kf\ d\mu+(\mathrm{H})\hspace{-0.4pc}\int_Kg\ d\mu$.
\item $(\mathrm{H})\hspace{-0.4pc}\int_{A\cup^*B}f\ d\mu=(\mathrm{H})\hspace{-0.4pc}\int_Af\ d\mu+(\mathrm{H})\hspace{-0.4pc}\int_Bf\ d\mu$ if $A,B\in{\cal S}$.
\item $(\mathrm{H})\hspace{-0.4pc}\int_Kf\ d\mu=(0,0)$ iff $f=(0,0)\ \mu$ almost everywhere.\qed
\end{enumerate}
\end{prp}

\begin{prp}\label{p3b}Let $(K,{\cal S},\mu)$ be an h-measure space. Let $f:K\to[0,+\infty)\times[0,+\infty)$ be a measurable simple function. Then for $L\in{\cal S}$ let $\nu(L)=(\mathrm{H})\hspace{-0.4pc}\int_Lf\ d\mu$. Then $\nu$ is an h-measure on ${\cal S}$.\qed
\end{prp}

As in subsection \cite{lahi2}\ref{hi2-ss1}, the next step is to define the integral of non-negative functions.

\begin{df}\label{d6}Let $(K,{\cal S},\mu)$ be an h-measure space. Let ${f:K\to[0,+\infty)\times[0,+\infty)}$ be a measurable function. Then set
\[(\mathrm{H})\hspace{-0.4pc}\int\limits_Kf\ d\mu=\sup\left\{(\mathrm{H})\hspace{-0.4pc}\int\limits_Kg\ d\mu:g\text{ is a measurable simple function}\text{ and }(0,0)\leq g\leq f\right\},\]
where a simple function takes its values from $[0,+\infty)\times[0,+\infty)$.\qed
\end{df}

\begin{prp}\label{p15}Let $(K,{\cal S},\mu)$ be an h-measure space. Let ${f,g:K\to[0,+\infty)\times[0,+\infty)}$ be measurable functions. Then the followings hold.
\vspace{-0.6pc}\begin{enumerate}\setlength\itemsep{-0.3em}
\item If $f$ is a simple function, then definition \ref{d6} provides the same result than definition \ref{d5}.
\item $(0,0)\leq(\mathrm{H})\hspace{-0.4pc}\int_Kf\ d\mu$.
\item If $f\leq g$, then $(\mathrm{H})\hspace{-0.4pc}\int_Kf\ d\mu\leq(\mathrm{H})\hspace{-0.4pc}\int_Kg\ d\mu$.
\item If $c\in[0,+\infty)\times[0,+\infty)$, then $(\mathrm{H})\hspace{-0.4pc}\int_Kcf\ d\mu=c\cdot(\mathrm{H})\hspace{-0.4pc}\int_Kf\ d\mu$.
\item $(\mathrm{H})\hspace{-0.4pc}\int_{A\cup^*B}f\ d\mu=(\mathrm{H})\hspace{-0.4pc}\int_Af\ d\mu+(\mathrm{H})\hspace{-0.4pc}\int_Bf\ d\mu$ if $A,B\in{\cal S},\ A\cap B=\emptyset$.
\end{enumerate}
\end{prp}
\begin{proof}(1),(2),(3),(4) can be proved as before in proposition \cite{lahi2}\ref{hi2-p13}.
\par (5): First let $\Phi_{f,A}={\{g:g\text{ is a measurable simple function},\ g|_{K-A}\equiv (0,0),\ (0,0)\leq g\leq f\}}$. Clearly $\Phi_{f,A}+\Phi_{f,B}\subset \Phi_{f,A\cup B}$ holds. Let $g\in\Phi_{f,A\cup B}$. Let $g_A(x)=g(x)$ if $x\in A$, and $g_A(x)=(0,0)$ if $x\notin A$. Define similarly $g_B$. Then evidently $g=g_A+g_B$ and $g_A\in\Phi_{f,A},\ g_B\in\Phi_{f,B}$ which yields that $\Phi_{f,A\cup B}\subset\Phi_{f,A}+\Phi_{f,B}$. Altogether $\Phi_{f,A}+\Phi_{f,B}=\Phi_{f,A\cup B}$ gives the claim.
\end{proof}

\begin{prp}\label{p16}Let $(K,{\cal S},\mu)$ be an h-measure space. Let ${f:K\to[0,+\infty)\times[0,+\infty)}$ be measurable function. Then 
$(\mathrm{H})\hspace{-0.4pc}\int_Kf\ d\mu=(0,0)$ iff $f=(0,0)\ \mu$ almost everywhere.
\end{prp}
\begin{proof}The sufficiency is obvious. 
\par To see the necessity, let $L=\{x\in K:\pi_1(f(x))>0\}$. Suppose that $\mu(L)>(0,0)$. Let $L_n=\left\{x\in L:\pi_1(f(x))>\frac{1}{n}\right\}$. Clearly $L=\bigcup_{n=1}^{\infty}L_n$, hence there is $n\in\N$ such that $\mu(L_n)>(0,0)$. Let $g(x)=(0,0)$ if $x\notin L_n$, and $g(x)=\left(\frac{1}{n},0\right)$ if $x\in L_n$. Obviously $(0,0)\leq g\leq f$ and $g$ is a measurable simple function, thus $(0,0)<(\mathrm{H})\hspace{-0.4pc}\int_Kg\ d\mu\leq(\mathrm{H})\hspace{-0.4pc}\int_Kf\ d\mu$ which is a contradiction. Therefore $\mu(L)=(0,0)$.
\par Let $M=\{x\in K-L:\pi_2(f(x))>0\}$. Again suppose indirectly that $\mu(M)>(0,0)$. Let $M_n=\left\{x\in L:\pi_2(f(x))>\frac{1}{n}\right\}$. Clearly $M=\bigcup_{n=1}^{\infty}M_n$, hence there is $n\in\N$ such that $\mu(M_n)>(0,0)$. Let $g(x)=(0,0)$ if $x\notin M_n$, and $g(x)=\left(0,\frac{1}{n}\right)$ if $x\in M_n$. Obviously $(0,0)\leq g\leq f$ and $g$ is a measurable simple function, thus $(0,0)<(\mathrm{H})\hspace{-0.4pc}\int_Kg\ d\mu\leq(\mathrm{H})\hspace{-0.4pc}\int_Kf\ d\mu$ which is a contradiction. Therefore $\mu(M)=(0,0)$.
\par Altogether we got that $\mu(\{x\in K:f(x)>(0,0)\})=(0,0)$.
\end{proof}


Now we are going to prove that the integral is additive.

\begin{lem}\label{l6}Let $(K,{\cal S},\mu)$ be an h-measure space. Let ${f:K\to[0,+\infty)\times[0,+\infty)}$ be a measurable function such that for all $x\in K\ \pi_1(f(x))=d$. Let $\tilde{f}(x)=\pi_2(f(x))\ (x\in K)$. Then
\[(\mathrm{H})\hspace{-0.4pc}\int\limits_Kf\ d\mu=(d,0)+(\mathrm{H})\hspace{-0.4pc}\int\limits_K\tilde{f}\ d\mu.\]
\end{lem}
\begin{proof}Let ${g:K\to[0,+\infty)\times[0,+\infty)}$ be a measurable simple function such that $(0,0)\leq g\leq f$. Then $\pi_1(g)\leq d$. If $x\in K$, then let 
\[g'(x)=
\begin{cases}
g(x)&\text{if }\pi_1(g(x))=d\\
(d,0)&\text{if }\pi_1(g(x))<d.
\end{cases}\]
Then $g'$ is a measurable simple function such that $(0,0)\leq g\leq g'\leq f$. Let $\tilde{g}(x)=\pi_2(g'(x))\ (x\in K)$. Then $\tilde{g}:K\to[0,+\infty)$ is a measurable simple function such that $0\leq \tilde{g}\leq \tilde{f}$. Evident calculation shows for both $d=0$ and $d>0$ that
\[(\mathrm{H})\hspace{-0.4pc}\int\limits_Kg\ d\mu\leq(\mathrm{H})\hspace{-0.4pc}\int\limits_Kg'\ d\mu\leq(d,0)+(\mathrm{H})\hspace{-0.4pc}\int\limits_K\tilde{g}\ d\mu\leq (d,0)+(\mathrm{H})\hspace{-0.4pc}\int\limits_K\tilde{f}\ d\mu.\]
This gives that $(\mathrm{H})\hspace{-0.4pc}\int_Kf\ d\mu\leq (d,0)+(\mathrm{H})\hspace{-0.4pc}\int_K\tilde{f}\ d\mu$.
\par To see the opposite case, let $\tilde{g}:K\to[0,+\infty)$ is a measurable simple function such that $0\leq \tilde{g}\leq \tilde{f}$. Let $g:K\to[0,+\infty)\times[0,+\infty)\ g(x)=(d,\tilde{g}(x))\ (x\in K)$. Then $g$ is a measurable simple function such that $(0,0)\leq g\leq f$. Then evidently
\[(d,0)+(\mathrm{H})\hspace{-0.4pc}\int\limits_K\tilde{g}\ d\mu\leq(\mathrm{H})\hspace{-0.4pc}\int\limits_Kg\ d\mu\leq(\mathrm{H})\hspace{-0.4pc}\int\limits_Kf\ d\mu,\]
which yields that $(d,0)+(\mathrm{H})\hspace{-0.4pc}\int_K\tilde{f}\ d\mu\leq(\mathrm{H})\hspace{-0.4pc}\int_Kf\ d\mu$.
\end{proof}

\begin{lem}\label{l7}Let $(K,{\cal S},\mu)$ be an h-measure space. Let ${f,g:K\to[0,+\infty)\times[0,+\infty)}$ be measurable functions such that for all $x\in K\ {\pi_1(f(x))=\pi_1(g(x))=d}$. Then
\[(\mathrm{H})\hspace{-0.4pc}\int\limits_Kf+g\ d\mu=(\mathrm{H})\hspace{-0.4pc}\int\limits_Kf\ d\mu+(\mathrm{H})\hspace{-0.4pc}\int\limits_Kg\ d\mu.\]
\end{lem}
\begin{proof}By lemma \ref{l6}, we get that
\[(\mathrm{H})\hspace{-0.4pc}\int\limits_Kf+g\ d\mu=(d,0)+(\mathrm{H})\hspace{-0.4pc}\int\limits_K\widetilde{f+g}\ d\mu=(d,0)+(\mathrm{H})\hspace{-0.4pc}\int\limits_K\tilde{f}+\tilde{g}\ d\mu=\]
\[(d,0)+(\mathrm{H})\hspace{-0.4pc}\int\limits_K\tilde{f}\ d\mu+(\mathrm{H})\hspace{-0.4pc}\int\limits_K\tilde{g}\ d\mu=(d,0)+(\mathrm{H})\hspace{-0.4pc}\int\limits_K\tilde{f}\ d\mu+(d,0)+(\mathrm{H})\hspace{-0.4pc}\int\limits_K\tilde{g}\ d\mu=\]
\[(\mathrm{H})\hspace{-0.4pc}\int\limits_Kf\ d\mu+(\mathrm{H})\hspace{-0.4pc}\int\limits_Kg\ d\mu,\]
applying $\widetilde{f+g}=\tilde{f}+\tilde{g}$ and proposition \cite{lahi2}\ref{hi2-p17} as well.
\end{proof}

\begin{lem}\label{l5}Let $(K,{\cal S},\mu)$ be an h-measure space. Let ${f,g:K\to[0,+\infty)\times[0,+\infty)}$ be measurable functions. Then $\{x\in K:\pi_1(f(x))<\pi_1(g(x))\}\in{\cal S}$.
\end{lem}
\begin{proof}Obviously
\[\{x\in K:\pi_1(f(x))<\pi_1(g(x))\}=\bigcup_{r\in\mathbb{Q}}\big(\{x\in K:f(x)<(r,0)\}\cup\{x\in K:(r,+\infty)<g(x)\}\big).\qedhere\]
\end{proof}

\begin{thm}Let $(K,{\cal S},\mu)$ be an h-measure space. Let ${f,g:K\to[0,+\infty)\times[0,+\infty)}$ be measurable functions. Then 
\[(\mathrm{H})\hspace{-0.4pc}\int\limits_Kf+g\ d\mu=(\mathrm{H})\hspace{-0.4pc}\int\limits_Kf\ d\mu+(\mathrm{H})\hspace{-0.4pc}\int\limits_Kg\ d\mu.\]
\end{thm}
\begin{proof}If $(\mathrm{H})\hspace{-0.4pc}\int_Kf\ d\mu=(\mathrm{H})\hspace{-0.4pc}\int_Kg\ d\mu=(0,0)$, then we are done by \ref{p16}. Hence we can suppose that any of those integrals does not equal to $(0,0)$.
\par Let $\Phi_h=\{j:(0,0)\leq j\leq h,\ j\text{ is a measurable simple function}\}$. Then $\Phi_f+\Phi_g\subset\Phi_{f+g}$ gives that $(\mathrm{H})\hspace{-0.4pc}\int_Kf\ d\mu+(\mathrm{H})\hspace{-0.4pc}\int_Kg\ d\mu\leq(\mathrm{H})\hspace{-0.4pc}\int_Kf+g\ d\mu$.
\par Let $A_1=\{x\in K:\pi_1(f(x))<\pi_1(g(x))\}, {A_2=\{x\in K:\pi_1(f(x))>\pi_1(g(x))\}}$, ${A_3=\{x\in K:\pi_1(f(x))=\pi_1(g(x))\}}$. Clearly $K=A_1\cup^* A_2\cup^* A_3$ and they are measurable by lemma \ref{l5}. By proposition \ref{p15}(5) it is enough to show for $i=1,2,3$ that we have
\[(\mathrm{H})\hspace{-0.4pc}\int\limits_{A_i}f+g\ d\mu=(\mathrm{H})\hspace{-0.4pc}\int\limits_{A_i}f\ d\mu+(\mathrm{H})\hspace{-0.4pc}\int\limits_{A_i}g\ d\mu.\]
First we show that for $A_1$. If $x\in A_1$, then $f(x)+g(x)=g(x)$, hence $(\mathrm{H})\hspace{-0.4pc}\int_{A_1}f+g\ d\mu=(\mathrm{H})\hspace{-0.4pc}\int_{A_1}g\ d\mu$, which together with the first inequality (at the beginning of the proof) gives the claim.
\par For $A_2$ we get the statement in the same way.
\par For $A_3$, let $h$ be a measurable simple function on $A_3$ such that $h=\sum_{i=1}^n(d_i,m_i)\chi_{S_i}\leq f+g$ and $(0,0)<(d_i,m_i),\ (0,0)<\mu(S_i)$. If $x\in S_i$, then $d_i\leq \pi_1(f(x)+g(x))=\pi_1(f(x))=\pi_1(g(x))$. For $1\leq i\leq n$ let 
\[S_i^0=\{x\in S_i:d_i=\pi_1(f(x))\},\ S_i^1=\{x\in S_i:d_i<\pi_1(f(x))\}.\]
Let $S^0=\bigcup_{i=1}^nS_i^0,\ S^1=\bigcup_{i=1}^nS_i^1$. Clearly $A_3=S^0\cup^* S^1$. Again by proposition \ref{p15}(5) it is enough to show for $i=0,1$ that we have
\[(\mathrm{H})\hspace{-0.4pc}\int\limits_{S^i}f+g\ d\mu=(\mathrm{H})\hspace{-0.4pc}\int\limits_{S^i}f\ d\mu+(\mathrm{H})\hspace{-0.4pc}\int\limits_{S^i}g\ d\mu.\]
By lemma \ref{l7}, it is true for $S^0$.
\par To see it for $S^1$, we define simple functions $s_i^1,t_i^1$ on $S^1$ $(1\leq i\leq n)$. If $x\in S_i^1$, let $s_i^1(x)=t_i^1(x)=\left(d_i,\frac{m_i}{2}\right)$, and $s_i^1(x)=t_i^1(x)=(0,0)$ when $x\in S^1-S^1_i$.
\par Let $s=\sum_{i=1}^ns_i^1\chi_{S_i^1},\ t=\sum_{i=1}^nt_i^1\chi_{S_i^1}$. Then $s\in \Phi_{f,S^1}, t\in \Phi_{g,S^1}$ and $h|_{S^1}=s+t$. By proposition \ref{p6b}(4) it gives that
\[(\mathrm{H})\hspace{-0.4pc}\int\limits_{S^1}h\ d\mu=(\mathrm{H})\hspace{-0.4pc}\int\limits_{S^1}s\ d\mu+(\mathrm{H})\hspace{-0.4pc}\int\limits_{S^1}t\ d\mu\leq (\mathrm{H})\hspace{-0.4pc}\int\limits_{S^1}f\ d\mu+(\mathrm{H})\hspace{-0.4pc}\int\limits_{S^1}g\ d\mu\]
which yields that
\[(\mathrm{H})\hspace{-0.4pc}\int\limits_{S^1}f+g\ d\mu\leq (\mathrm{H})\hspace{-0.4pc}\int\limits_{S^1}f\ d\mu+(\mathrm{H})\hspace{-0.4pc}\int\limits_{S^1}g\ d\mu\]
which together with the first inequality (at the beginning of the proof) gives the claim.
\end{proof}

\begin{thm}Let $(K,{\cal S},\mu)$ be an h-measure space. Let $f:K\to[0,+\infty)\times[0,+\infty)$ be a measurable function. Then 
\[\nu(L)=(\mathrm{H})\hspace{-0.4pc}\int\limits_Lf\ d\mu\ \ (L\in{\cal S})\]
is an h-measure on ${\cal S}$.
\end{thm}
\begin{proof}Let $(K_n)$ be a pairwise disjoint sequence of sets from ${\cal S}$. We can assume that $K=\cup_{n=1}^{\infty}K_n$.
\par For $n\in\N$ let $f_n(x)=f(x)$ if $x\in K_1\cup\dots\cup K_n$, otherwise let $f_n(x)=(0,0)$. Clearly $(0,0)\leq f_n\leq f$. By proposition \ref{p15}(3) and (5), we get that
\[(\mathrm{H})\hspace{-0.4pc}\int\limits_Kf\ d\mu\geq\sum_{n=1}^{\infty}(\mathrm{H})\hspace{-0.4pc}\int\limits_{K_n}f\ d\mu.\]
Let $(0,0)\leq g\leq f$ be a measurable simple function. Evidently $g|_{K_n}$ is a measurable simple function, hence $(\mathrm{H})\hspace{-0.4pc}\int_{K_n}g\ d\mu\leq(\mathrm{H})\hspace{-0.4pc}\int_{K_n}f\ d\mu$. By proposition \ref{p3b}, we get that
\[(\mathrm{H})\hspace{-0.4pc}\int\limits_Kg\ d\mu=\sum_{n=1}^{\infty}(\mathrm{H})\hspace{-0.4pc}\int\limits_{K_n}g\ d\mu\leq\sum_{n=1}^{\infty}(\mathrm{H})\hspace{-0.4pc}\int\limits_{K_n}f\ d\mu\]
which yields that $(\mathrm{H})\hspace{-0.4pc}\int_Kf\ d\mu\leq\sum_{n=1}^{\infty}(\mathrm{H})\hspace{-0.4pc}\int_{K_n}f\ d\mu$.
\end{proof}

\begin{ex}It is not true that for every measurable function $f$, there is a sequence $(g_n)$ of measurable simple functions such that $g_n\leq g_{n+1}\ (n\in\N)$ and $g_n\to f$.
\par Let $K=(0,1),\ {\cal S}$ be the Borel sets of $K$.
\par Let $f:K\to[0,1]\times[0,1],\ f(x)=(x,x)$. Clearly $f$ is measurable. Suppose that $f$ can be approximated by simple functions in the above mentioned way. Let ${H=\{\pi_1(\mathrm{Ran\ }g_n):n\in\N\}}\subset[0,1]$. Clearly $H$ is countable, hence there is $x\in K-H$. Then there is no $n\in\N$ such that $(x,0)<g_n(x)<(x,1)$, but that interval is a neighborhood of $f(x)=(x,x)$.\qed
\end{ex}

\begin{ex}It is not true that if $f_n\to f,\ (0,0)\leq f_n\leq f_{n+1}\ (n\in\N)$, then $(\mathrm{H})\hspace{-0.4pc}\int_Kf_n\ d\mu\to(\mathrm{H})\hspace{-0.4pc}\int_Kf\ d\mu$.
\par Let $K=(0,1),\ {\cal S}$ be the Borel sets of $K$, and let $\mu=(1,\lambda)$ where $\lambda$ is the Lebesgue measure on $K$. By proposition \cite{lahi2}\ref{hi2-p12}, $(K,{\cal S},\mu)$ is an h-measure space. 
\par Let $f_n:K\to[0,1]\times[0,1],\ f_n(x)=(\sqrt[n]{x},\sqrt[n]{x})\ (n\in\N, x\in K)$. Clearly $f_n$ is measurable $(n\in\N)$. We show that $(\mathrm{H})\hspace{-0.4pc}\int_Kf_n\ d\mu=(2,0)\ (n\in\N)$. Let $s=\sum_{i=1}^n(d_i,m_i)\chi_{A_i}$ be a measurable simple function such that $(0,0)\leq s\leq f_n$. Then $\sum_{i=1}^n(d_i,m_i)\mu(A_i)\leq\sum_{i=1}^n(1,0)(1,1)\leq(2,0)$, therefore $(\mathrm{H})\hspace{-0.4pc}\int_Ks\ d\mu\leq(2,0)$ and $(\mathrm{H})\hspace{-0.4pc}\int_Kf_n\ d\mu\leq(2,0)$.
\par Now for $\varepsilon>0$, let $s(x)=(1-\varepsilon,0)$ if $\sqrt[n]{x}>1-\varepsilon$, otherwise let $s(x)=(0,0)$. Clearly $s$ is a measurable simple function such that $(0,0)\leq s\leq f_n$. We get that $(\mathrm{H})\hspace{-0.4pc}\int_Ks\ d\mu=(2-\varepsilon,0)$, which now gives that $(\mathrm{H})\hspace{-0.4pc}\int_Kf_n\ d\mu=(2,0)$.
\par But $f(x)=(1,1)\ (x\in K)$, thus $(\mathrm{H})\hspace{-0.4pc}\int_Kf\ d\mu=(2,1)$.\qed
\end{ex}

\begin{rem}\label{r1}By remark \ref{r2}, this integral is a natural generalization of the one described in \cite{lahi2}.
\par Furthermore it also handles the case when $f$ taking values from ${[0,+\infty)\times[0,+\infty)}$ and $\mu$ is an ordinary measure only: just apply proposition \cite{lahi2}\ref{hi2-p12}. The most natural embedding is when we have $\mu(S)=(0,\nu(S))\ (S\in{\cal S})$ for given measure $\nu$ on ${\cal S}$.
\end{rem}

Now we present a way how the integral can be calculated. First we need an extension of the notion of simple function where we allow that its second coordinate function takes $+\infty$ as a value as well.

\begin{df}We call $g:K\to[0,+\infty)\times[0,+\infty]$ an i-simple function if $\mathrm{Ran}g$ is finite.
\end{df}

\begin{lem}Let $(K,{\cal S},\mu)$ be an h-measure space. Let ${f:K\to[0,+\infty)\times[0,+\infty)}$ be a measurable function. Then 
\[(\mathrm{H})\hspace{-0.4pc}\int\limits_Kf\ d\mu=\sup\left\{(\mathrm{H})\hspace{-0.4pc}\int\limits_Kg\ d\mu:g\text{ is a measurable i-simple function}\text{ and }(0,0)\leq g\leq f\right\}.\]
\end{lem}
\begin{proof}Let $s$ denote the right hand side. 
\par As the set for which the $\sup$ regards is greater, we get that ${(\mathrm{H})\hspace{-0.4pc}\int_Kf\ d\mu\leq s}$.
\par Let $g:K\to[0,+\infty)\times[0,+\infty]$ be an i-simple function, $g=\sum_{i=1}^n(d'_i,m'_i)\chi_{A_i}$ such that $(0,0)<(d'_i,m'_i)\in[0,+\infty)\times[0,+\infty]$ and $A_i\cap A_j=\emptyset\ (i\ne j)$. Let $\mu(A_i)=(d''_i,m''_i)>(0,0)$. For $k\in\N$ and $i\in\{1,\dots,n\}$ let
\[n_{ik}=\begin{cases}
m'_i&\text{if }m'_i<+\infty\\
1&\text{if }m'_i=+\infty\text{ and }m''_i=0\\
k&\text{if }m'_i=+\infty\text{ and }m''_i>0.
\end{cases}\]
Then let $g_k=\sum_{i=1}^n(d'_i,n_{ik})\chi_{A_i}$. Clearly 
\[(\mathrm{H})\hspace{-0.4pc}\int\limits_{A_i}g_k\ d\mu=(d'_i+d''_i,n_{ik}m''_i)\to(d'_i+d''_i,m'_im''_i)=(\mathrm{H})\hspace{-0.4pc}\int\limits_{A_i}g\ d\mu\ \ \ (i\in\{1,\dots,n\}).\]
As the addition in $[0,+\infty)\times[0,+\infty]$ is continuous (\cite{lahi}2.5), we get that
\[(\mathrm{H})\hspace{-0.4pc}\int\limits_Kg_k\ d\mu\to(\mathrm{H})\hspace{-0.4pc}\int\limits_Kg\ d\mu\]
which gives that $s\leq(\mathrm{H})\hspace{-0.4pc}\int_Kf\ d\mu$.
\end{proof}

\begin{thm}\label{t4}Let $(K,{\cal S},\mu)$ be an h-measure space. Let ${f:K\to[0,+\infty)\times[0,+\infty)}$ be a measurable function. Then
\[(\mathrm{H})\hspace{-0.4pc}\int\limits_Kf\ d\mu=(d,m)\]
where 
\[d=\sup\left\{d'+d'':\exists L\subset K,\ \mu(L)=(d'',m'')>(0,0),\ \inf_L f=(d',m')>(0,0)\right\},\]
\vspace{-0.8pc}\begin{multline*}m=\sup\bigg\{\sum_{i=1}^nm'_im''_i:\exists L_1,\dots,L_n\subset K,\ L_i\cap L_j=\emptyset\ (i\ne j),\ \mu(L_i)=(d''_i,m''_i)>(0,0),\\
\inf_{L_i} f=(d'_i,m'_i)>(0,0),\ d=d'_i+d''_i\ \ (1\leq i\leq n)\bigg\},
\end{multline*}
with the convention that $\sup\emptyset=0$.
\end{thm}
\begin{proof}Let $(d,m)=(\mathrm{H})\hspace{-0.4pc}\int_Kf\ d\mu$, 
\newline$\hat{d}=\sup\left\{d'+d'':\exists L\subset K,\ \mu(L)=(d'',m'')>(0,0),\ \inf_L f=(d',m')>(0,0)\right\}$, 
\newline$\hat{m}=\sup\{\sum_{i=1}^nm'_im''_i:\exists L_1,\dots,L_n\subset K,\ L_i\cap L_j=\emptyset\ (i\ne j),\ \mu(L_i)=(d''_i,m''_i)>(0,0),\ \inf_{L_i} f=(d'_i,m'_i)>(0,0),\ d=d'_i+d''_i\ \ (1\leq i\leq n)\}$.
\par If $L\subset K, \inf_L f=(d',m')>(0,0)$, then let $g(x)=(d',m')$ if $x\in L$, and $g(x)=(0,0)$ if $x\in K-L$. Let $\mu(L)=(d'',m'')>(0,0)$. Then $g$ is a measurable i-simple function such that $(0,0)\leq g\leq f$, hence
\[d'+d''=\pi_1\left((\mathrm{H})\hspace{-0.4pc}\int\limits_Kg\ d\mu\right)\leq\pi_1\left((\mathrm{H})\hspace{-0.4pc}\int\limits_Kf\ d\mu\right)=d\]
which gives that $\hat{d}\leq d$.
\par For $\varepsilon>0$ let $g:K\to\R$ be a measurable simple function such that $(0,0)\leq g\leq f$, $g=\sum_{i=1}^n(d'_i,m'_i)\chi_{A_i}$ such that $(d'_i,m'_i)>(0,0)$ and $A_i\cap A_j=\emptyset\ (i\ne j)$ and 
\[\pi_1\left((\mathrm{H})\hspace{-0.4pc}\int\limits_Kf\ d\mu\right)-\varepsilon<\pi_1\left((\mathrm{H})\hspace{-0.4pc}\int\limits_Kg\ d\mu\right)=\max\{d'_i+d''_i:1\leq i\leq n\}\]
where $\mu(A_i)=(d''_i,m''_i)>(0,0)$.
Therefore there is $i\in\{1,\dots,n\}$ such that $d-\varepsilon<d'_i+d''_i$. For that $i$ we have that $(0,0)<(d'_i,m'_i)\leq\inf_{A_i} f$ which gives that $d-\varepsilon<\hat{d}$. As it holds for all $\varepsilon$, we get that $d\leq \hat{d}$.
\smallskip
\par When validating the equality for the second coordinate, there are two cases.
\vspace{-0.3pc}\begin{enumerate}\setlength\itemsep{-0.3em}
\item If $L\subset K, \inf_L f=(d',m')>(0,0), \mu(L)=(d'',m'')$, then $d'+d''<d$. Then it can be readily seen that $m=\hat{m}=0$.
\item Suppose that there are $L\subset K, \inf_L f=(d',m')>(0,0), \mu(L)=(d'',m'')$ such that $d'+d''=d$. Then let $L_1,\dots,L_n\subset K$ such that $L_i\cap L_j=\emptyset\ (i\ne j),\ \mu(L_i)=(d''_i,m''_i)>(0,0),\ \inf_{L_i} f=(d'_i,m'_i)>(0,0)$ and for all $i\in\{1,\dots,n\}$ $d=d'_i+d''_i$. Let $g(x)=(d'_i,m'_i)$ if $x\in L_i$, and $g(x)=(0,0)$ if $x\in K-\cup_iL_i$. Then $g$ is a measurable i-simple function such that $(0,0)\leq g\leq f$, hence
\[\left(d,\sum_{i=1}^nm'_im''_i\right)=(\mathrm{H})\hspace{-0.4pc}\int\limits_Kg\ d\mu\leq(\mathrm{H})\hspace{-0.4pc}\int\limits_Kf\ d\mu=(d,m)\]
which gives that $\hat{m}\leq m$.
\par For $\varepsilon>0$ let $g:K\to\R$ be a measurable simple function such that $(0,0)\leq g\leq f$, $g=\sum_{i=1}^n(d'_i,m'_i)\chi_{A_i}$ such that $(d'_i,m'_i)>(0,0)$ and $A_i\cap A_j=\emptyset\ (i\ne j)$ and $\mu(A_i)=(d''_i,m''_i)>(0,0)$ and for all $i\in\{1,\dots,n\}$ $d=d'_i+d''_i$ and 
\[\pi_2\left((\mathrm{H})\hspace{-0.4pc}\int\limits_Kf\ d\mu\right)-\varepsilon<\pi_2\left((\mathrm{H})\hspace{-0.4pc}\int\limits_Kg\ d\mu\right)=\sum_{i=1}^nm'_im''_i\leq\hat{m}.\]
Hence $m-\varepsilon<\hat{m}$. As it holds for all $\varepsilon$, we get that $m\leq \hat{m}$.\qedhere 
\end{enumerate}
\end{proof}

\begin{rem}Note that in the formulation of the theorem, inside the definition of $d$, the condition $\mu(L)>(0,0)$ is essential, but inside the definition of $m$, the condition $\mu(L_i)>(0,0)$ can be omitted.
\end{rem}

We can apply the previous theorem when $\mu$ is an ordinary measure.

\begin{lem}\label{l2}Let $(K,{\cal S},\nu)$ be a measure space. Let $f:K\to[0,+\infty)$ be a measurable function. Then 
\[\mathrm{ess\ sup\ }f=\sup\{\inf_L f:L\in{\cal S},\ \nu(L)>0\}.\]
\end{lem}
\begin{proof}Let $s=\mathrm{ess\ sup\ }f$ and $s'=\sup\{\inf_L f:L\in{\cal S},\ \nu(L)>0\}$. We know that $s=\mathrm{ess\ sup\ }f=\inf\{a\geq 0: \nu(f^{-1}([a,\infty)))=0\}$. 
\par If $L\in{\cal S},\ \nu(L)>0$, then let $a=\inf_L f$. Then $L\subset f^{-1}([a,\infty))$ which yields that $\nu(f^{-1}([a,\infty)))>0$ hence $a\leq s$ which gives that $s'\leq s$.
\par Let $a<s$. Then $\nu(f^{-1}([a,\infty)))>0$. Let $L=f^{-1}([a,\infty)),\ b=\inf_L f$. But as $s'\geq b\geq a$ holds for all such $a$, we get that $s'\geq s$.
\end{proof}

\begin{thm}Let $(K,{\cal S},\nu)$ be a measure space. Let ${f:K\to[0,+\infty)\times[0,+\infty)}$ be a measurable function. Then
\[(\mathrm{H})\hspace{-0.4pc}\int\limits_Kf\ d\mu=\left(\mathrm{ess\ sup\ }\pi_1(f),\int\limits_{E_f}\pi_2(f)\ d\nu\right)\]
where $E_f=\{x\in K:\pi_1(f(x))=\mathrm{ess\ sup\ }\pi_1(f)\}$.
\end{thm}
\begin{proof}Let $\mu=(0,\nu)$ be the associated h-measure (see also remark \ref{r1}). 
\par Applying theorem \ref{t4} (also using the same notation), we get that 
\newline$d=\sup\left\{d':\exists L\subset K,\ \nu(L)>0,\ \inf_L f=(d',m')>(0,0)\right\}$ which equals to $\mathrm{ess\ sup\ }\pi_1(f)$ by lemma \ref{l2} and using that $\pi_1(\inf_L f)=\inf_L \pi_1(f)$.
\par For $m$ we get that
\vspace{-0.8pc}\begin{multline*}m=\sup\bigg\{\sum_{i=1}^nm'_im''_i:\exists L_1,\dots,L_n\subset K,\ L_i\cap L_j=\emptyset\ (i\ne j),\ \nu(L_i)=m''_i>0,\\
\inf_{L_i} f=(d,m'_i)>(0,0)\ (1\leq i\leq n)\bigg\}.
\end{multline*}
Let $L'_i=\{x\in L_i:\pi_1(f(x))=d\}, L''_i=\{x\in L_i:\pi_1(f(x))>d\}$. Let $L''_{in}=\left\{x\in L_i:\pi_1(f(x))>d+\frac{1}{n}\right\}\ (n\in\N)$. By the above result on $d$, we get that $\nu(L''_{in})=0$ for all $n$. Then $L''_i=\bigcup_{n=1}^{\infty}L''_{in}$ gives that $\nu(L''_i)=0$. Hence we can replace $\inf_{L_i} f=(d,m'_i)$ with $\pi_1(f|_{L_i})=d,\ \inf_{L_i}\pi_2(f)=m'_i$. But that simply means that $L_i\subset E_f$ and $m=\int_{E_f}\pi_2(f)\ d\nu$.
\end{proof}

\section{Applications}

We provide three applications that are sophisticated deficiency measurements for continuity, lineness and convexity (for basic notions and notations see \cite{lamd} and \cite{lahi} subsection 3.1). Here we will just define the deficiency measurements with a little explanation, and we will not add reasonable statements as our purpose is not that currently.

\smallskip

Let $f:\mathbb{R}\to\mathbb{R}$ be a function. We want to measure how much $f$ is not continuous. Set
\[\mathrm{defi}(f,\text{''continuity''},\text{''H-int cluster''})=(\mathrm{H})\hspace{-0.4pc}\int\limits_{\mathbb{R}}\dm\big(\Lambda_f(x)-\{f(x)\}\big)\ d\dm,\]
where $\Lambda_f(x)=\big\{y\in\mathbb{R}:\exists (x_n)\text{ sequence such that } x_n\to x, f(x_n)\to y\big\}$.
That is we calculate the ''Hausdorff-size'' of each cluster set $\Lambda_f(x)$, i.e. its Hausdorff-dimension and its associated measure, and then we ''summarize'' them using the Hausdorff-integral. In \cite{lamd} we show that 
this deficiency significantly differs from $\dm\Big(\cup\big\{\{x\}\times\Lambda_f(x):x\in\mathbb{R}\big\}-\mathrm{graph }f\Big)$.
That example also shows that the Fubini theorem does not hold for the Hausdorff-integral in general.


\smallskip

Let $K\subset\mathbb{R}^2$ be a measurable set. We want to measure how much $K$ is not a line. If $y$ is a point of the line $e$, then $f(y)$ denotes the line perpendicular to $e$ at $y$, and $t(y)=K\cap f(y)-\{y\}$. Set
\[\mathrm{defi}(f,\text{''lineness''},\text{''H-int out''})=\inf\left\{(\mathrm{H})\hspace{-0.4pc}\int\limits_e\dm(t(y))\ d\dm:e\text{ is a line in }\mathbb{R}^2\right\}\]

\smallskip

Let $K\subset\mathbb{R}^2$ be bounded, measurable. We want to measure how much $K$ is not convex. Set
\[\mathrm{defi}(K,\text{''convex''},\text{''mhd''})=(\mathrm{H})\hspace{-0.4pc}\int\limits_{K\times K}\dm(\overline{xy}-K)\ d\dm(x,y),\]
where $\overline{xy}$ is the line segment between $x$ and $y$, and the domain of $\dm$ is the subsets of $\mathbb{R}^4$.


\medskip

{\footnotesize

\noindent

\noindent email: alosonczi1@gmail.com\\
}
\end{document}